\documentclass[12pt]{article}
\usepackage{eucal,amssymb,amsmath,graphicx, amsfonts, latexsym,epsfig,makeidx,subfigure,multirow,hhline}
\usepackage[T1]{fontenc}

\begin{document} \parskip=5pt plus1pt minus1pt \parindent=0pt
\title{A weighted configuration model and\linebreak inhomogeneous epidemics}
\author{Tom Britton\thanks{Department of Mathematics, Stockholm University, 106 91 Stockholm, Sweden; {\tt tomb at math.su.se, mia at math.su.se}} \and Maria Deijfen$^*$\and Fredrik Liljeros\thanks{Department of Sociology, Stockholm University, 106 91 Stockholm, Sweden; {\tt liljeros at sociology.su.se}}}
\date{April 2011}
\maketitle

\begin{abstract}
\noindent A random graph model with prescribed degree distribution and degree dependent edge weights is introduced. Each vertex is independently equipped with a random number of half-edges and each half-edge is assigned an integer valued weight according to a distribution that is allowed to depend on the degree of its vertex. Half-edges with the same weight are then paired randomly to create edges. An expression for the threshold for the appearance of a giant component in the resulting graph is derived using results on multi-type branching processes. The same technique also gives an expression for the basic reproduction number for an epidemic on the graph where the probability that a certain edge is used for transmission is a function of the edge weight. It is demonstrated that, if vertices with large degree tend to have large (small) weights on their edges and if the transmission probability increases with the edge weight, then it is easier (harder) for the epidemic to take off compared to a randomized epidemic with the same degree and weight distribution. A recipe for calculating the probability of a large outbreak in the epidemic and the size of such an outbreak is also given. Finally, the model is fitted to three empirical weighted networks of importance for the spread of contagious diseases and it is shown that $R_0$ can be substantially over- or underestimated if the correlation between degree and weight is not taken into account.

\noindent
\vspace{0.3cm}

\vspace{0.2cm}


\end{abstract}

\section{Introduction}

Random graph models aimed at describing large scale network structures have been much studied the last few years; see e.g.\ \cite{Dur,vdH} and the references therein. One reason for being interested in such models is that realistic network models make it possible to quantify and predict the outcome of epidemics taking place on the networks. Most work so far on network epidemics have been restricted to unweighted graphs where transmission takes place along all edges with the same probability; see e.g.\ \cite{A-98,A-99,BJM07, Neal}. The purpose of the present work is to formulate a model for a weighted network where the weight of an edge is a function of the degrees of the adjacent vertices. Moreover, a simple epidemic on such a network will be analyzed where the transmission probability is taken to be a function of the edge weight. This generalizes work by Britton et al.\ \cite{BJM07} and Deijfen \cite{D}.

The network model we shall work with is a generalization of the well known configuration model; see \cite{MR-95,MR-98}. Each vertex is independently assigned a number of half-edges according to a fixed degree distribution and each half-edge is given an integer valued weight according to a distribution that is allowed to depend on the degree of its vertex. Half-edges with the same weight are then paired completely at random to create edges. Once the graph has been generated we let a Reed-Frost epidemic with weight-dependent infection probabilities spread on it. More specifically, each vertex that is infective at time $t$ ($t=1,2,\ldots$) independently infects each of its susceptible neighbors with a probability that is a function $\pi(w)$ of the weight $w$ of the connecting edge and is then removed from the epidemic process. At time $t+1$, the infected neighbors become infective and transmit the infection according to the same rules. The weights hence have a twofold effect: they affect the structure of the underlying network and they determine the infection probabilities in the epidemic taking place on the network.

We analyse the model as the vertex population size $n$ tends to infinity using branching process approximations. A \emph{large outbreak} in the epidemic is said to occur if a positive fraction of the population is asymptotically infected and the \emph{basic reproduction number}, $R_0$, is defined as a function of the parameters of the model such that a large outbreak has positive probability if and only if $R_0>1$. Clearly, positive probability of a large outbreak in a Reed-Frost network epidemic initiated by one single infective is equivalent to the existence of a giant component -- that is, a component of order $n$ -- in a thinned version of the underlying network where each edge is kept independently with a probability specified by the infection probability through the edge. We derive an expression for the basic reproduction number in our model. Taking the infection probability deterministically equal to 1 in this expression yields a threshold parameter for the occurrence of a giant component in the underlying graph. Furthermore, we briefly mention how an expression for the probability of a large outbreak and the relative size of the outbreak can be derived. Indeed, these quantities coincide and can be characterized using the same branching process approximation that yields the basic reproduction number.

The basic reproduction number for our model is a function of the degree distribution, the (degree dependent) weight distribution and the (weight dependent) infection probability. We illustrate its behavior in a few specific cases. For instance we demonstrate that, in a case where the edges of high degree vertices tend to be given small weights (that is, the network has a negative degree-weight correlation) and the infection probability increases with the weight, it is more difficult for an epidemic to take off compared to a situation where the weights are assigned independently of the degrees, with the same fraction of edges of each weight. Conversely, if the edges of high degree vertices tend to have large weights (that is, the network has positive degree-weight correlation), epidemics take off more easily. We also illustrate our results on three empirical networks. The first network comprises patients that have been registered at the same ward at a hospital and the weights here correspond to the number of days that pairs of patients have spent simultaneously at the same ward. The second network is obtained from register data on household and workplace structure in Sweden. Vertices represent workplaces and the weight on a link between two workplaces indicates the number of households that contain one person working at one of the workplaces and one person working at the other. The third data network consist of sexual contacts and the weights here represent the number of intercourses per contact. The first two networks have positive degree-weight correlation whereas the third data set has negative degree-weight correlation.

The rest of the paper is organized as follows. In Section 2, we define the network model, comment on its relation to the standard configuration model and define the epidemic model on the network. Section 3 is devoted to the derivation of the basic reproduction number for the epidemic and, as a special case, the threshold parameter for the appearance of a giant component in the network itself. The behavior of the threshold is then illustrated in a few examples in Section 3.1. In Section 4, we briefly describe how the probability of a large outbreak and, equivalently, the size of the outbreak, can be characterized. Section 5 contains applications with data from real-world networks, and Section 6 concludes with a short discussion on further work.

\section{Description of the model}

In this section we define the construction of the network model and the epidemic model, starting with the network.

\subsection{A model for a weighted network}\label{wei-net-mod}

Let $n$ denote the number of vertices. The vertices are first given strictly positive i.i.d.\ degrees $D_1,\dots ,D_n$ according to a prescribed distribution $P(D=d)=p(d)$ (the degrees are assumed strictly positive since vertices with degree zero are not part of the network and can therefore be disregarded). We think of $D_i$ as the number of stubs (or half-edges) sticking out of vertex $i$ ($i=1,\dots ,n$). Each stub is then independently assigned a non-negative integer valued weight, where the weights of the stubs of vertex $i$ have a distribution that is determined by the degree $D_i$ of the vertex. More specifically, if $D_i=d$, then the $d$ stubs are given i.i.d.\ weights $W_{i1}, \dots ,W_{id}$ according to a distribution $\{q(w|d)\}$, that is, we have
$$
P(W_{ij}=w|D_i=d)=q(w|d),\ j=1,\dots ,d.
$$
The network is now formed by paring up stubs with the same weight completely at random. More specifically, to pair up the stubs with weight $w$, first pick two stubs at random from the set of all stubs with weight $w$ and join them into an edge, then pick two stubs at random from the set of remaining stubs with weight $w$ and join them. And so on. If the number of stubs with weight $w$ is odd, we throw away the last stub. This pairing procedure is applied separately for each weight $w$.

To conclude, the model is defined in a very similar way as the well-known configuration model using the degree distribution $\{p(d)\}$ and the family of weight distributions $\{ q(w|d)\}$. Indeed, we retrieve the standard configuration model if all stubs are deterministically assigned the same weight regardless of the degree. If, more generally, the edges are assigned weights independently of the degrees -- that is, if $q(w|d)=q(w)$ for all $d$ -- then the model is equivalent to first generating the graph according to the standard configuration model and then assigning the weights independently according to the distribution $q(w)$. This case has previously been analyzed in \cite{D}.

It is not hard to see that the fact that we are removing the last remaining stub with a given weight in case the total number of stubs with this weight is odd does not affect the degree and weight distributions in the network in the limit as $n\to\infty$. Furthermore, for the standard configuration model, it is well-known that the fraction of self-loops and multiple edges between vertices is small as soon as the degree distribution has finite second moment. More specifically, removing self-loops and multiple edges does not change the degree distribution in the graph, and the probability of obtaining a simple graph is bounded away from 0 as $n\to\infty$; see \cite[Theorem 7.9]{vdH} and \cite[Lemma 5.5]{BJM07}. Clearly the same applies to our generalized model. Furthermore, the graph is tree-like, meaning that with high probability it does not contain short cycles. This allows for various types of branching process approximations; see e.g.\ \cite{A-98,BJM07,JL} for rigorous treatments of the standard configuration model.

Before proceeding, we remark that the proposed model can give rise to graphs that are very different from the standard configuration model. For instance, if all stubs at vertices with degree $d$ are assigned weight $w_d$, with $w_{d_1}\neq w_{d_2}$ for $d_1\neq d_2$, then a graph is obtained where edges run only between vertices with the same degree and where the set of vertices with degree $d$ constitute a sub-graph with the same structure as a graph obtained from the standard configuration model with constant degree $d$. The configuration model yields a fully connected graph (in probability) if $P(D\geq 3)=1$, see \cite[Theorem 10.14]{vdH}, and hence the above assignment of weights leads to a graph where, for each $d=3,4,\ldots$, all vertices with degree $d$ constitute a connected component of their own. This contrasts with the standard configuration model where a giant component -- that is, a component containing a positive fraction of the vertices -- is unique and occurs if and only if $E[D(D-1)]/E[D]>1$; see e.g.\ \cite{JL,MR-95}. Various generalizations of the above division of the graph based on the degrees of the vertices are possible. For instance, take subsets $A_1,A_2,\ldots\subset \mathbb{N}$ and $B_1,B_2\ldots\subset\mathbb{N}$, and assign the weights in such a way that vertices with degrees in $A_i$ are given weights in $B_i$. If $A_i\cap A_j=B_i\cap B_j=\emptyset$ for $i\neq j$, then each component in the graph will contain only vertices with degrees in a given set $A_i$.

\subsection{The epidemic model}

Given the network defined above, a Reed-Frost type epidemic spreads randomly in the network in the following way. Initially, at time $t=1$, one randomly selected vertex is infected and the remaining vertices are susceptible. The infection then spreads in generations in that a vertex that is infected at time $t$ ($t=1,2,\ldots$) infects each of its susceptible neighbors independently with a probability that depends on the weight of the connecting edge: if the weight is $w$, then the infection probability is $\pi (w)\in[0,1]$. At time $t+1$, all vertices that were infected at time $t$ become immune (or die) and play no further role in the epidemic. The epidemic goes on until there are no new infections -- then the epidemic stops, and the vertices that have been infected during the course of the outbreak make up the set of ultimately infected vertices.

As for the infection probability $\pi (w)$, there are many possible choices. Typically, we think of the weight of an edge as representing the strength or the intensity of the connection, and the infection probability $\pi(w)$ is then increasing in $w$. For instance, in a situation where the integer weight $w$ on an edge represents the number of contacts between the corresponding vertices in a certain time interval, one natural candidate is
\[
\pi(w)=1-(1-s)^w,
 \]
where $s$ is the per-contact probability of infection. However it is also possible to let the weights represent resistances in the connections, and $\pi(w)$ is then naturally decreasing. Finally, if we let $\pi (w)\equiv 1$ for all $w$ there is no thinning and we retrieve the original network.

\section{Threshold parameter}\label{sec:R0}

In this section we derive an asymptotic expression for the basic reproduction number in the epidemic as $n\to\infty$ or, equivalently, for the threshold for the appearance of a giant component in the weighted network where an edge with weight $w$ is removed with probability $1-\pi(w)$. The expression is valid under the assumption that a component can contain vertices with all possible degrees. This excludes cases where the weights are assigned in such a way that certain degrees are isolated in separate components, as described in Section 2.2. Furthermore, in order to be able to rely on standard results for multi-type branching processes, we shall assume throughout that the degree distribution as well as the weight distributions have bounded support. We expect however that the resulting expressions are valid as soon as the degree distribution has finite second moment.

Denote a vertex with degree $d$ by `$d$-vertex' and let $p_d(k)$ denote the probability that a given neighbor of a $d$-vertex is a $k$-vertex. A $d$-vertex has $d$ neighbors and, because the weights of the stubs are independent, the degrees of different neighbors become independent as $n\to\infty$. The probability that a given neighbor has degree $k$ depends on the weight $w$ of the connecting edge: if the weight is $w$, the probability that it leads to a $k$-vertex equals the probability $\tilde{p}_w(k)$ that a randomly chosen $w$-stub belongs to a $k$-vertex, where we have
$$
\tilde{p}_w(k)=\frac{q(w|k)kp(k)}{\sum_j q(w|j)jp(j)}.
$$
Hence
\begin{equation}
p_d(k)=\sum_wq(w|d)\tilde{p}_w(k).\label{pdk}
\end{equation}

Now consider the epidemic process and assume that $n$ is large. Except for the index case, all infected vertices will have been infected by one of their neighbors. On the other hand, because there are with high probability no short cycles in the network, the remaining neighbors of an infected vertex during the early stages of the epidemic will be susceptible with high probability. As a consequence, an infected $d$-vertex will with high probability have $d-1$ susceptible neighbors during the early stages. How many of these will the $d$-vertex infect and what are the degrees of the infected neighbors? The answer depends on the degree distribution of the neighbors and on the weights of the connecting edges. For $k=1,2,\ldots$, let $p^x_d(k)$ denote the probability that a given susceptible neighbor has degree $k$ and becomes infected by our $d$-vertex and write $X_d(k)$ for the number of susceptible $k$-neighbors that become infected. As in deriving (\ref{pdk}), we obtain
\begin{equation}
p_d^x(k)= \sum_w\pi(w)q(w|d)\tilde{p}_w(k),\label{pxdk}
\end{equation}
and conclude that $E[X_d(k)]=(d-1)p_d^x(k)$.

During the early stages of the epidemic, when the fraction of already infected vertices is still negligible, the epidemic may be approximated by a multitype branching process, where the type of a vertex is given by its degree. For vertices in generation $t\geq 2$, the expected offspring is specified above (the offspring of the index case is different; see Section 4). In particular, during the early stages of the epidemic (excluding the first generation), the expected number of $k$-vertices that an infected $d$-vertex infects equals
\[
m_{dk}=E[X_d(k)]=(d-1)p_d^x(k)=(d-1)\sum_w\pi (w)q(w|d)\tilde{p}_w(k).
\]
The matrix $M=(m_{dk})_{d,k\geq 1}$ is known as the mean offspring matrix and it is well-known from the theory of multi-type branching processes that, under the assumption that each type has the possibility of giving rise to offspring of any other type within a finite number of generations, the process has a positive probability of growing beyond all limits if and only if the largest eigenvalue of the offspring matrix exceeds 1; see e.g.\ \cite[Chapter 4]{J}. In our setting, the type 1 vertices are infertile and the above assumption is hence not met. It is however easy to see that our process has a positive probability of exploding if and only if the process defined by the type $d$ vertices, with $d\geq 2$, has a positive probability of exploding. Indeed, the degree 1 vertices do not contribute to further spread of the epidemic (they do however contribute to the final size of the epidemic; see Section 4). The basic reproduction number $R_0$, is hence given by the largest eigenvalue of the matrix $M_2=(m_{dk})_{d,k\geq 2}$ and a major outbreak can occur if and only if $R_0>1$. Taking $\pi(w)\equiv 1$ gives the threshold for the occurrence of a giant component in the original network.

Before proceeding, we remark that $M$ can be written as the product $A\cdot B$ of two matrices $A=(a_{dw})$ and $B=(b_{wk})$, where $a_{dw}=(d-1)\pi(w)q(w|d)$ is the expected number of edges with weight $w$ from a $d$-vertex that are used to transmit infection, and $b_{wk}=\tilde{p}_w(k)$ is the probability that a given edge with weight $w$ is attached to a $k$-vertex. The exploration of the degrees of the neighbors that are infected by a given $d$-vertex can hence be divided in two steps: first the $d$-vertex gives rise to a number of weighted transmission links, and then the other end of each transmission link is connected to a vertex whose degree depends on the weight of the link.

\subsection{Examples}

Here we investigate the threshold parameter for some specific choices of degree and weight distributions.

\noindent \textbf{Example 3.1} First we consider the case where the weights are independent of the degrees, that is, $q(w|d)=q(w)$ for all $d$. As already pointed out, the model is then equivalent to first generating the graph according to the standard configuration model and then assigning a weight independently to each edge according to the weight distribution $q(w)$, a case previously analyzed in \cite{D}. Using (\ref{pdk}), we get that
\begin{equation}\label{sb}
p_d(k)=\frac{kp(k)}{\sum_jjp(j)}=\frac{kp(k)}{\mu_D}:=\tilde{p}(k),
\end{equation}
so the degree distribution of the neighbors of a $d$-individual is independent of $d$ and given by the size-biased degree distribution. The mean offspring matrix $M=(m_{dk})$ is given by
$$
m_{dk}=(d-1)\frac{kp(k)}{\mu_D}\sum_w\pi(w)q(w)=(d-1)\tilde{p}(k)E[\pi(W)],
$$
where $E[\pi(W)]$ is the unconditional transmission probability of a randomly selected edge (which is independent of the adjacent vertices in this example). Since this matrix can be written as a column vector multiplied by a row vector, the largest eigenvalue equals the sum of the diagonal elements, that is,
$$
R_0=E[\pi(W)]\sum_k(k-1)\tilde{p}(k) = E[\pi(W)]\left( \mu_D+\frac{\sigma^2_D-\mu_D}{\mu_D}\right),
$$
where $\mu_D$ and $\sigma_D^2$ denote the mean and variance, respectively, in the degree distribution. The threshold is hence the same as for a homogeneous infection on the standard configuration model with infection probability given by the expected transmission probability of the connections, which is in agreement with the result in \cite{D}. Furthermore, when $W\equiv 1$, the transmission probability along any edge is $\pi (1)=:\pi$ and the above expression conincides with previous results for homogeneous epidemics on the configuration model; see e.g.\ \cite{BJM07}.\hfill$\Box$\medskip

\noindent \textbf{Example 3.2} When the weights are correlated with the degrees, the underlying graph can have a structure that is very different from graphs obtained from the standard configuration model. Consider for instance a case where the only possible degrees are 1 and 3 and with $p(1)=1-p(3)=0.8$ (and $\pi(w)\equiv 1$). The value of the critical parameter $E[D(D-1)]/E[D]$ for the standard configuration model is then 0.86. On the other hand, assume that there are two possible weights $w_1$ and $w_2$ for the 1-vertices and two possible weights $w_2$ and $w_3$ for the 3-vertices. Figure 1 shows a plot of the critical parameter as a function of $q(w_2|3)$ for $q(w_2|1)=0.3$ (so $q(w_1|1)=0.7$ and $q(w_3|3)=1-q(w_2|3)$). We see that the graph is supercritical for all choices of $q(w_2|3)$. This is because 1-vertices and 3-vertices are now only connected by edges with weight $w_2$ and, when the fraction of edges with weight $w_2$ on the 1-vertices is small, the interference from the 1-vertices is not enough to suppress the giant component formed by the 3-vertices without the presence of 1-vertices. Note that, since $\pi(w)\equiv 1$ in this example, the only effect of the weight is that they introduce degree correlation in the graph.\hfill$\Box$\medskip

\begin{figure}[p]
\begin{center}
\mbox{\epsfig{file=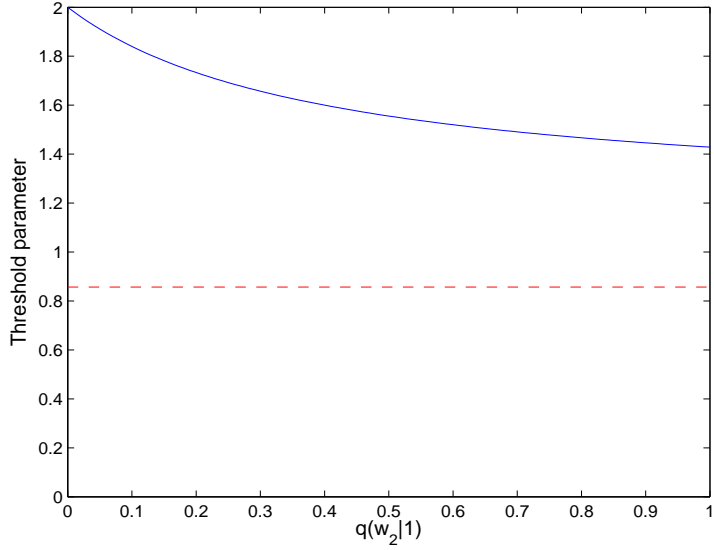,width=11cm, height=8cm}}
\end{center}
\caption{Threshold parameter for $p(1)=1-p(3)=0.8$ with two possible weights $w_1$ and $w_2$ for degree 1 and weights $w_2$ and $w_3$ for degree 3, and with $\pi(w)\equiv 1$ (cf.\ Example 3.2). The value of the threshold is plotted as a function of $q(w_2|3)$ for $q(w_2|1)=0.3$ (solid line) and the horizontal line indicates the threshold for the standard configuration model with the same degree distribution. }
\end{figure}

\noindent \textbf{Example 3.3} We next give an example with Po(4)-degree distribution conditioned on being in [1,200], and with two possible weights $w_1$ and $w_2$, with $w_1<w_2$. First take $q(w_2|d)=d^{-\alpha}$ (and $q(w_1|d)=1-q(w_2|d)$) for $\alpha>0$-- that is, the probability that a given link from a $d$-individual has the larger weight $w_2$ decays with $d$ at rate $\alpha$ -- and set $\pi(w_1)=0.1$ and $\pi(w_2)=0.7$. Figure 2 shows a plot of $R_0$ (solid line) against $\alpha$. The dashed line represents the basic reproduction number for an epidemic where the weights are assigned independently of the degrees in such a way that the fraction of edges with weight $w_1$ and $w_2$ respectively is the same as in the above network, that is, the probability of assigning weight $w_2$ to a given stub is set to
$$
q(w_2)=\sum_k\tilde{p}(k)k^{-\alpha},
$$
where $\tilde p(k)=kp(k)/\mu_D$.
The plot reveals that the epidemic with negative degree-correlated weights has a smaller $R_0$, which is explained by the fact that the high-weight edges are then less likely to be connected to high-degree vertices. Naturally, the reproduction numbers converge as $\alpha$ increases to the threshold for a homogeneous epidemic with infection probability 0.1 (since there will be very few high-weight edges for large $\alpha$).

Figure 3 shows a plot of $R_0$ against $\alpha$ for the same setup, but with $q(w_2|d)=1-d^{-\alpha}$, that is, the probability that a given link from a $d$-individual has the larger weight $w_2$ now instead increases with $d$. Now the epidemic with the degree-correlated weights has the largest $R_0$, since the large-weight edges occur with higher probability at high-degree vertices, making it easier for the epidemic to take off.\hfill$\Box$\medskip

\begin{figure}[p]
\begin{center}
\mbox{\epsfig{file=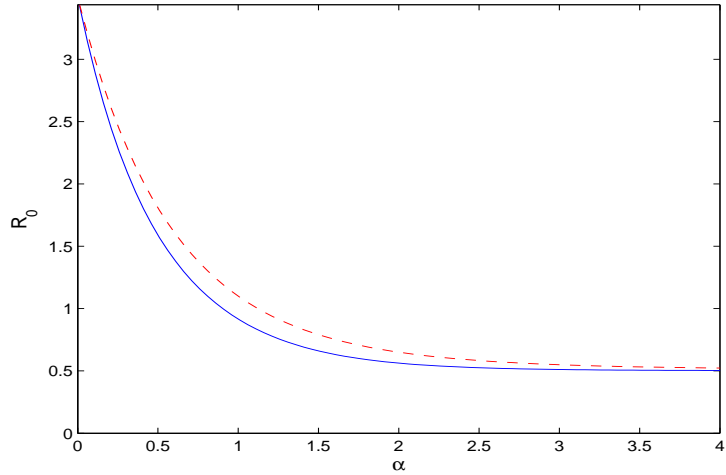,width=11cm, height=7cm}}
\end{center}
\caption{Basic reproduction numbers for two-point weights $w_1<w_2$ with $\pi(w_1)=0.1$ and $\pi(w_2)=0.7$ and with Po(4) degrees conditioned on being in [1,200] (cf.\ Example 3.3). Degree correlated weights with $P(w_2|D=d)=d^{-\alpha}$ (solid line) and degree independent weights with the same fraction of large-weight edges (dashed line) plotted against $\alpha$.}
\end{figure}

\begin{figure}[p]
\begin{center}
\mbox{\epsfig{file=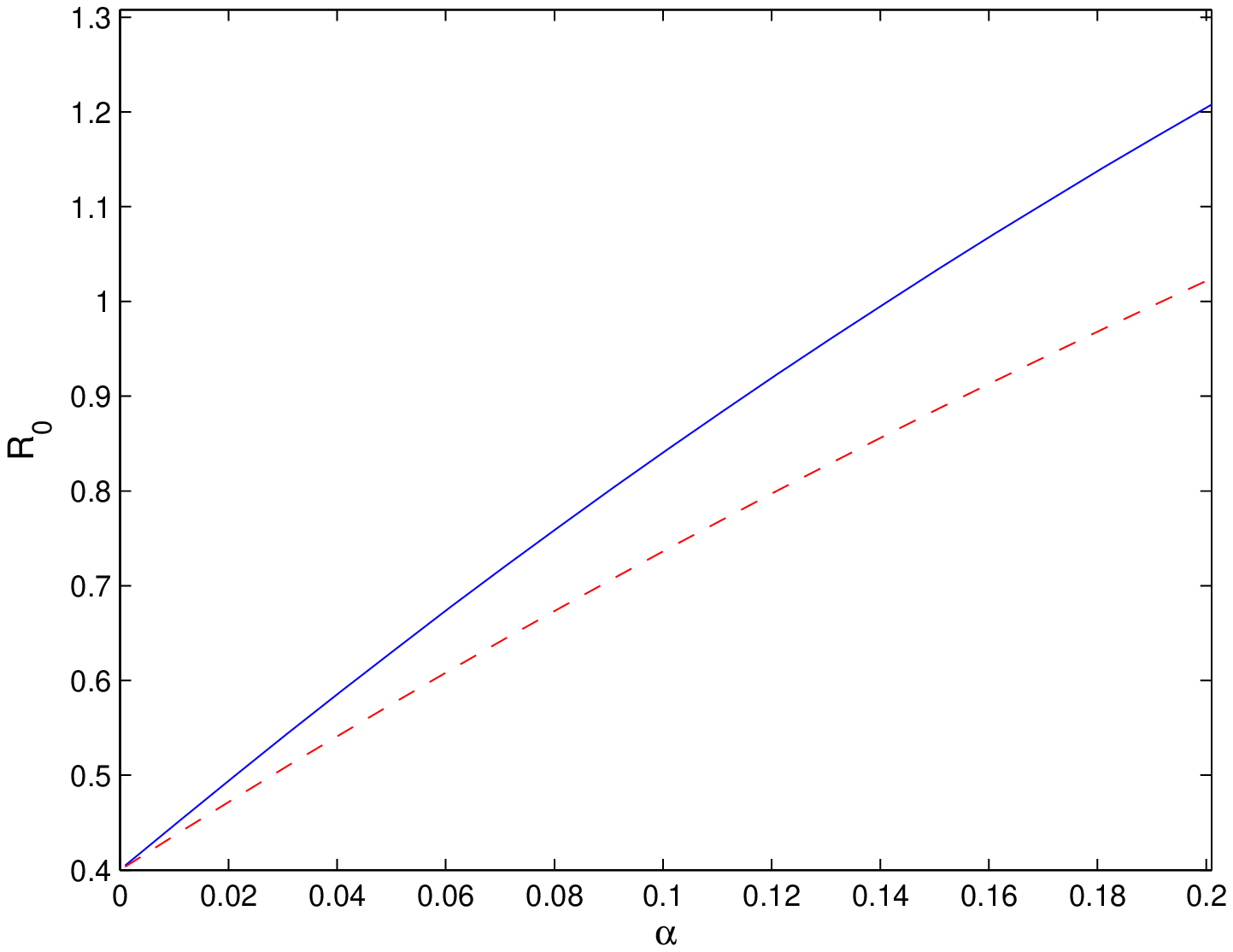,width=11cm, height=8cm}}
\end{center}
\caption{Basic reproduction number plotted against $\alpha$ for Example 3.3, but with $P(w_2|D=d)=1-d^{-\alpha}$.}
\end{figure}

\noindent \textbf{Example 3.4} Finally consider a case where the expected total weight of all edges of a vertex, given the degree of the vertex, is fixed and independent of the degree. More specifically, fix $\gamma\in\mathbb{R}$ and, for $w\in\mathbb{N}$, let $q(w|d)$ be the probability that a Po($\gamma/d$)-variable takes on the value $w$. The expected weight on an edge of a $d$-vertex is then $\gamma/d$ and hence the expected total weight of a $d$-vertex is $\gamma$. The transmission probability for an edge with weight $w$ is set to $\pi(w)=1-(1-s)^{w+1}$. Figure 4 shows a plot of the basic reproduction number against $s$ when the degree distribution is Po(8), conditioned on being positive, and with $\gamma=8$ (so on average 8 friends and expected total weight $\gamma =8$). For comparison, a plot of the basic reproduction number for an epidemic where the weights are assigned independently of the degrees is also included. There, $q(w|d)$ is taken to be the probability that a Po($\gamma/\mu$) variable takes on the value $w$, which means that the expected total weight per vertex still equals $\gamma=8$ (and $E(D)=\mu=8$). The plot reveals however that the basic reproduction number for the epidemic with degree dependent weights is much smaller (in fact always below the critical value 1). This follows from the fact that the expected total weight of a vertex is kept fixed independently of the degree which implies that the edges of high degree vertices will be assigned smaller weights.\hfill$\Box$\medskip

\begin{figure}[p]
\begin{center}
\mbox{\epsfig{file=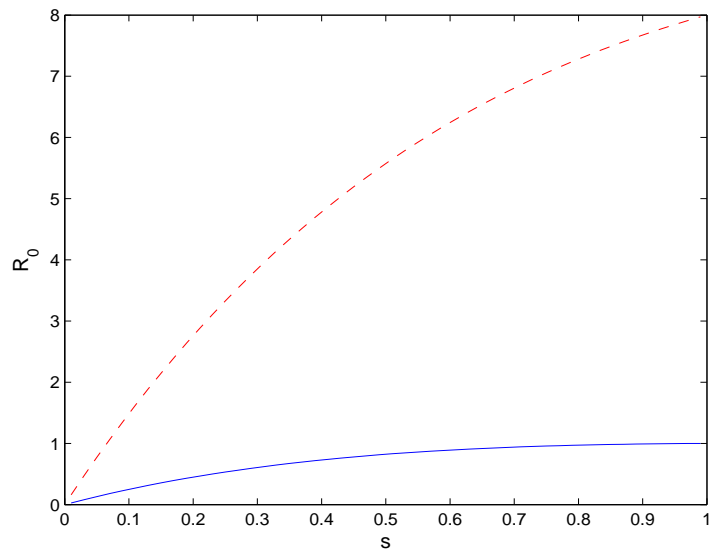,width=11cm, height=8cm}}
\end{center}
\caption{Plot of the basic reproduction number $R_0$ as a function of $s$ for  Example 3.4. The solid line represents the original example and dashed line the case when the weights are assigned independently of the degrees.}
\end{figure}

In Example 3.4, the expected total weight of a vertex conditionally on its degree, is kept fixed, which may be natural in many applications. From the perspective of the epidemic, a related, and perhaps even more important quantity, is the expected total ''infection pressure'' of a vertex. This is quantified for a $d$-vertex by
\begin{equation}
t(d):=(d-1)\sum_wq(w|d)\pi(w)=(d-1)E[\pi(W)|D=d],\label{td}
\end{equation}
that is, $t(d)$ is the expected number of neighbors that a $d$-vertex infects during the early stages of the epidemic. If $t(d)$ is increasing (decreasing) in $d$, vertices with high degree tend to cause more (fewer) new cases. Furthermore, it is easy to see that, if $t(d)=t$ for all degrees $d$, then $R_0=t$. Equation (\ref{td}) also illustrates that two ''competing'' factors determine whether a large outbreak is possible or not. One factor is the degree distribution, further emphasized by the size-biasing for infected vertices. This factor is also present in un-weighted networks, as described in Example 3.1, and it is well-known that a heavy tailed degree distribution gives rise to a large $R_0$. The other factor is the expected transmission probability as a function of the degree, or more generally, the \emph{distribution} of the transmission probability as a function of the degree. If $E[\pi (W)|D=d]$ also increases with $d$ this will make $R_0$ even larger, whereas the perhaps more likely scenario that $E[\pi (W)|D=d]$ decreases with $d$ will typically downplay the role of a heavy tailed degree distribution (cf.\ \cite{BNL} for an illustration of this phenomenon only allowing two different weights). If the transmission probability is very small for vertices with high degree, a heavy tailed degree distribution may even lead to a smaller $R_0$.

\section{Outbreak probability and final size}

As mentioned in the introduction there is a close relationship between
the \emph{final size} of a major outbreak and the \emph{probability} of a major
outbreak for Reed-Frost type epidemics. In particular, if the epidemic is started
by one randomly selected index case, then, as the population size
$n$ tends to infinity, the probability $\rho$ of a major outbreak coincides with the proportion $\tau$ that is infected in case a major outbreak occurs; see e.g.\ \cite{B10}. Below we outline how to derive $\rho$ having the above dual interpretation. We do this by taking the methodology presented in Section 3 one step further.

Recall that $X_d(k)$ is the (random) number of $k$-vertices that a
$d$-vertex infects during the early stages of the epidemic
(excluding the first generation since the index case has all its
neighbors susceptible instead of all but one). The
$d$-vertex has $d-1$ susceptible neighbors (it was infected by
one neighbor but the others are, with high probability, still susceptible) and of
these it might infect some $k$-neighbors, some $j$-neighbors
and so on. It might however also \emph{fail} to infect some
$k$-neighbors, some $j$-neighbors and so on. What happens along
each outgoing stub from the $d$-vertex is independent and
identically distributed, and the probability that a $k$-individual is
infected equals $p_d^x(k)$ defined in (\ref{pxdk}). If we let $Y$ count
the number of neighbors that do not get infected (irrespective of type)
it hence follows that the vector $\mathbf{Z}_d=(X_d(1), \dots , X_d(k_{\max}), Y)$
is multinomially distributed with parameters $d-1$ and $p^x_d(1), \dots
, p_d^x(k_{\max}),\ \bar p$, where $\bar p=1-\sum_kp^x_d(k)$. This holds
for $d=1,\dots ,k_{\max}$. These so-called offspring distributions $\mathbf{Z}_d,\ d=1, \dots ,
k_{\max}$, defines the multi-type branching process with which the initial phase of the epidemic may be
approximated by (for large $n$).

We now use results from branching process theory (see e.g.\ \cite{J}) to derive
the probability $\rho$ that the branching process grows beyond all limits
(which corresponds to a major outbreak in the epidemic). First we note
that $1$-vertices that are infected in the second or
succeeding generations have no effect on $\rho$ because they will never spread the infection onwards and hence have no effect on the probability of growing beyond all limits. We hence neglect
$1$-vertices in the branching process. Define for $d=2, \dots , k_{\max}$
\[f_d(s_2, \dots, s_{k_{\max}})= E\left(s_2^{X_d(2)}\cdot \dots \cdot
s_{k_{\max}}^{X_d(k_{\max})}\right)=\left(p^x_d(1)+\sum_{k\ge
    2}s_kp^x_d(k)+\bar p \right)^{d-1},
\]
where the last equality follows from properties of the multinomial
distribution ($f_d$ is more or less the probability generating
function of $\mathbf{Z}_d$). The extinction probability $q_d$ if starting with one
$d$-vertex in generation 2 or later ($d=2, \dots , k_{\max}$) is
then given by the smallest solution to the equation system
\begin{align}
q_2&=f_2(q_2,\dots ,q_{k_{\max}})\nonumber\\
& \vdots \label{q-vec}\\
q_{k_{\max}}&= f_{k_{\max}}(q_2,\dots ,q_{k_{\max}}).\nonumber
\end{align}

From branching process theory it is well-known that $q_k=1$ for all
$k$ if $R_0\le 1$, whereas if $R_0>1$ there is a unique solution
$\mathbf{q}=(q_2,\dots ,q_{k_{\max}})$ satisfying $0<q_k<0$ for all $k$. These
are the extinction probabilities of ancestors in the second generation and later. The distribution of the number of vertices infected in the second generation (i.e.\ by the index case) is different for two reasons: the index case has degree distribution $\{p(d)\}$ and can infect all (rather than all but one) of its neighbors, implying that also $1$-vertices can start a major outbreak. We obtain the outbreak probability $\rho$ by conditioning first on the degree of the index case and secondly on the number of infected vertices of each type
that the index case infects (that is, those that make up the second generation). For $d=1,\dots ,d_{k_{max}}$, let $(\tilde X_d(2), \dots ,\tilde X_d(k_{\max}))$ have the same distribution as  $(X_d(2), \dots , X_d(k_{\max}))$ with the only difference that $d-1$ is replaced by $d$ in the multinomial distribution. The probability $1-\rho$ of not having a major outbreak is then given by
\begin{equation}
1-\rho = \sum_{d\ge 1}p(d)E\left( q_2^{\tilde X_d(2)}\cdot \dots \cdot q_{k_{\max}}^{\tilde X_d(k_{\max})}\right) = \sum_{d\ge 1}\left(p^x_d(1)+\sum_{k\ge
    2}q_kp^x_d(k)+\bar p \right)^{d}.\label{rho-def}
\end{equation}
The numerical difficulties in deriving $\rho$ are two-fold. First $p^x_d(k)$ defined in (\ref{pxdk}) have to be computed for all $d$ and $k$, and then the equation system (\ref{q-vec}) defining the vector $\mathbf{q}$ has to be solved numerically. Once this is done, computing $\rho$ is straightforward using (\ref{rho-def}).

\section{Fitting the network model to data}

Suppose data is collected from a weighted network with $n$ vertices. The information is summarized by a table $\{ n(w,d);w=1, 2, \dots, d=1,2, \dots\}$, where $n(w,d)$ denotes the number of edges with weight $w$ that are connected to $d$-vertices. Furthermore, let $n(d)$ denote the number of $d$-vertices. We then have $\sum_wn(w,d)=dn(d)$, or equivalently $n(d)=\sum_wn(w,d)/d$, satisfying $\sum_dn(d)=n$.

The data give rise to natural estimates of the degree distribution $p(d)$ and the weight distribution $q(w|d)$:
\begin{align}
\hat p(d)&=\frac{n(d)}{n},\\
\hat q(w|d)&=\frac{n(w,d)}{\sum_vn(v,d)}.
\end{align}

We now fit the model to three empirical weighted networks with importance for the spread of contagious diseases.

\noindent \textbf{Network 1: Inpatients}. The first network consists of hospital inpatients connected to each other by weighted links, where the weight on a link between two patients indicates the number of days that these patients have been registered simultaneously at the same ward in the region of Stockholm \cite{LGH}. A patient leaving a ward the same day as another patient enters the ward induces a link with weight 1, a patient leaving the ward the day after another patient enters the ward induces a link with weight 2, and so on. The total sampling period is seven days. The network has a positive degree-weight correlation.

\noindent \textbf{Network 2: Workplaces}. The second network is based on register data on workplaces and households in Sweden \cite{Chen1,Chen2}. Vertices represent workplaces and the weight on a link between two workplaces indicates the number or households where one person works at one of the workplaces and one person works at the other. Here, a household is defined as a married couple (with or without children), or a unmarried couple having children together, that live in the same flat or house. Unmarried couples without children and other individuals sharing a household are not registered as households in the database and are therefore not taken into account in the analysis. The network exhibits positive degree-weight correlation.

\noindent \textbf{Network 3: Sex contacts}. The third network describes survey data from a detailed sexual study carried out on the island of Gotland in Sweden in 1988 \cite{BNL,NL}. The survey contained questions that make it possible to estimate both the number of partners and the number of intercourses per partner for each respondent (during a year). We will not take gender into account in the analysis since, in its current form, the model suggested here is not applicable to bipartite network data. The network has negative degree-weight correlation.

\begin{table}\centering
\begin{tabular}{|l|c|c|c|}
  \hhline{~|---}
  \multicolumn{1}{c|}{} & Inpatients & Workplaces & Sex contacts\\
  \hhline{-|---}
  $n$ & 8227 & 318362 & 477\\
  \cline{1-4}
   $\#$ Links & 133762 & 957778 & 1159\\
  \cline{1-4}
  Mean($d$) & 50.2 & 104.0 & 2.4\\
  \cline{1-4}
  Stdev($d$) & 31.7 & 305.0 & 2.1\\
  \cline{1-4}
  Deg.corr. & 0.29 & 0.10 & -0.34\\
  \cline{1-4}
  Mean($w$) & 4.4 & 1.1 & 40.7\\
  \cline{1-4}
  Stdev($w$) & 2.3 & 0.9 & 59353\\
  \hline
\end{tabular}
 \caption{Summary statistics for the empirical networks.}
\end{table}

Basic statistical information about the networks can be found in Table 1. All three networks display non scale-free degree distributions; see Figure 5. A difference however is that the patient and workplace networks have an assortative structure (that is, positive degree correlation) while the sexual contact network has a disassortative structure. Furthermore, the workplace network stands out from the other two in that it has a mean weight close to 1 (meaning that most connected workplaces are only connected through one common household).

\begin{figure}[p]
\begin{center}
\mbox{\epsfig{file=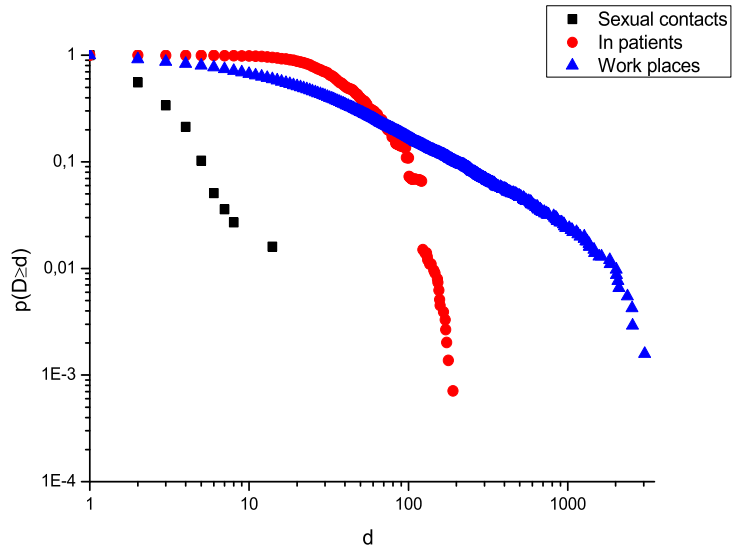,width=11cm, height=8cm}}
\end{center}
\caption{Degree tails for the empirical networks.}
\end{figure}

Figures 6-8 show $R_0$ as a function of $s$, where $\pi(w)=1-(1-s)^w$. In the three networks, the parameter $s$ is interpreted as risk of transmission per day spent at the same ward, risk of transmission per household connecting two workplaces and risk of transmission per sexual intercourse, respectively. As a reference, the values obtained in networks with the same degrees and weights but where the weights have been re-shuffled and assigned independently of the degree are also included.

\begin{figure}[p]
\begin{center}
\mbox{\epsfig{file=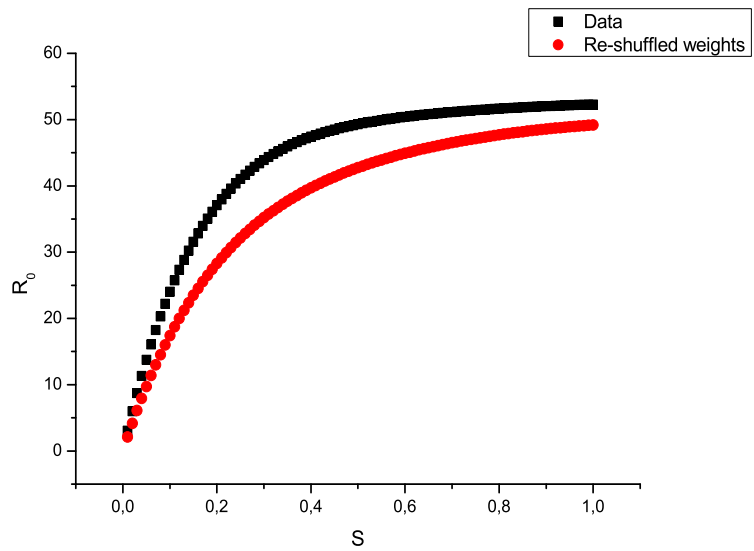,width=11cm, height=8cm}}
\end{center}
\caption{Basic reproduction number plotted against transmission probability per contact $s$ for the inpatient network.}
\end{figure}

\begin{figure}[p]
\begin{center}
\mbox{\epsfig{file=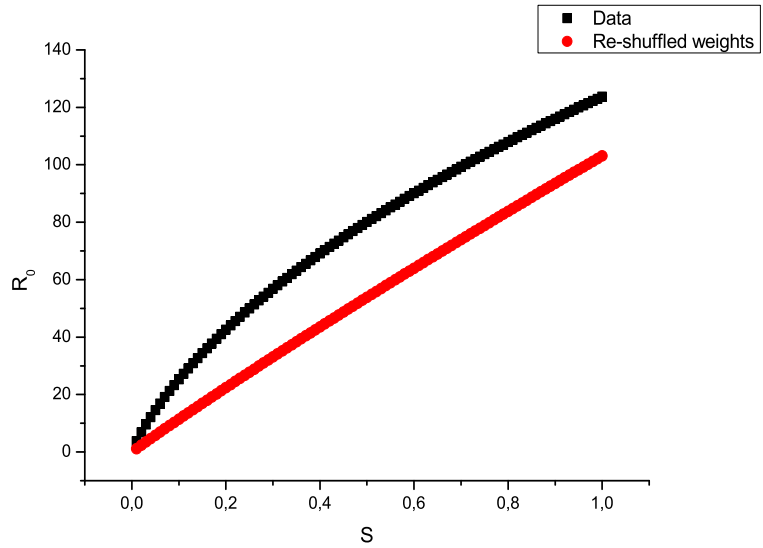,width=11cm, height=8cm}}
\end{center}
\caption{Basic reproduction number plotted against transmission probability per contact for the workplace network.}
\end{figure}

\begin{figure}[p]
\begin{center}
\mbox{\epsfig{file=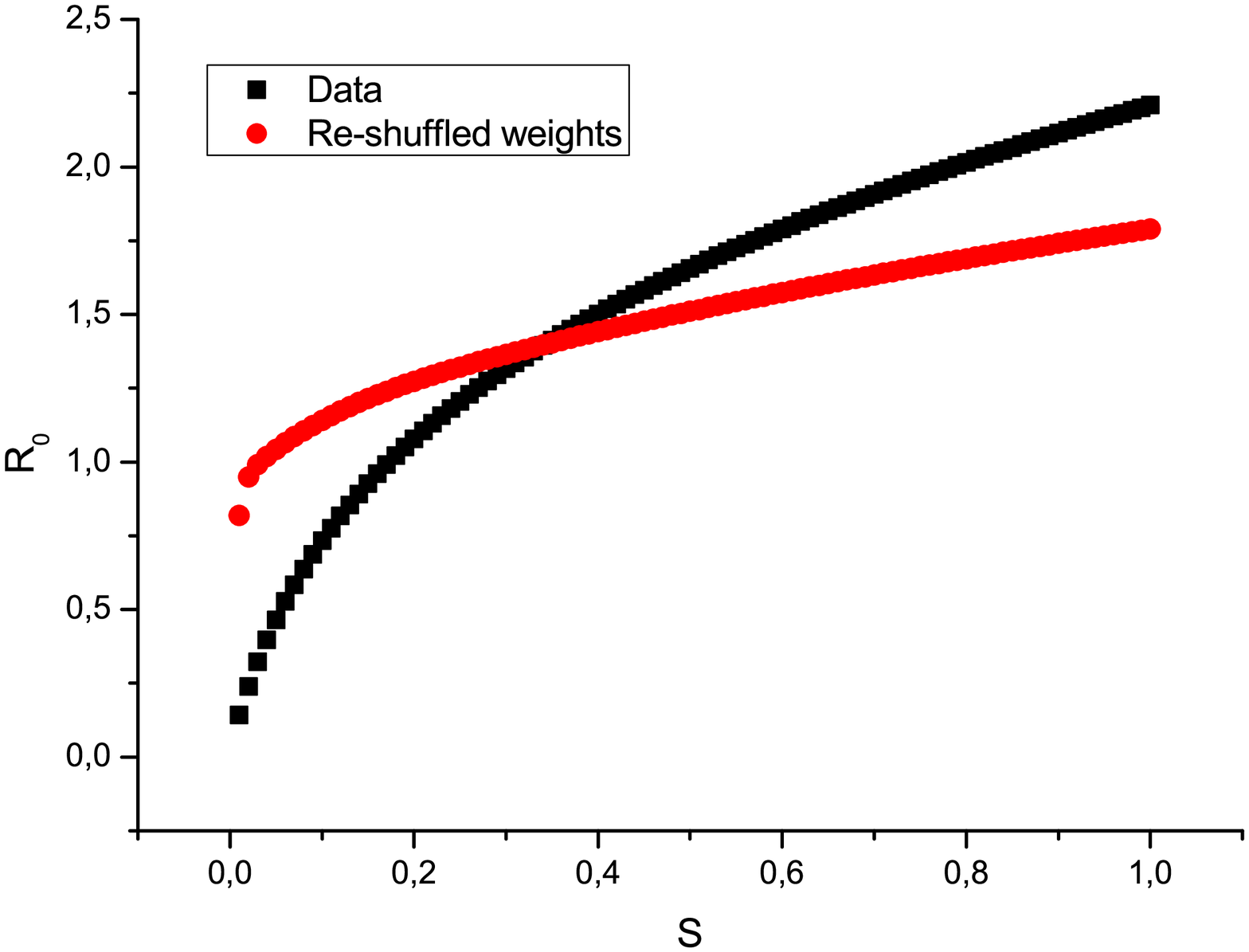,width=11cm, height=8cm}}
\end{center}
\caption{Basic reproduction number plotted against transmission probability per contact for the sexual contact network.}
\end{figure}

For both the hospital and the workplace networks we observe higher values of $R_0$ compared to the re-shuffled networks. This is because in those networks vertices with large degree have larger average weight on their edges compared to vertices with low degree, and will therefore be more effective spreaders than their counterparts in the re-shuffled networks. The sexual contact data exhibits a qualitatively different relationship between the empirical network and the re-shuffled network: here $R_0$ is smaller for the empirical data than the re-shuffled data when $s$ is small and larger than the re-shuffled data when $s$ is large. The reason for this is that large degree vertices will tend to have larger weights on their edges in the re-shuffled network, which is important when the risk of transmission per contact (that is, $s$) is small. On the other hand, the empirical data is more assortative than the re-shuffled data, which results in larger values of $R_0$ compared to the re-shuffled data when the risk of transmission per contact is high. Finally, the observed differences in functional form of $R_0$ as a function of $s$ between the networks can mainly be explained by the differences in average weight.

\section{Discussion}

In the present paper we have defined a network model with weighted edges and studied limiting properties of such a network as well as properties of an epidemic spreading on the network. The weights may have two effects: they may induce degree correlation in the underlying network, and they (typically) affect the transmission probability and hence the propagation of the epidemic. We have seen that the behavior of $R_0$ (as well as the outbreak probability) is affected by the joint distribution of the pair $(D, W)$, in particular, the correlation between $D$ and $W$ is important. The effects caused by the distribution of $(D,W)$ can then be either amplified or mitigated by the functional relationship between the weight and the transmission probability. The main conclusion is that $R_0$ can be substantially overestimated or underestimated if the degree-weight correlation is not taken into account.

Having relaxed the assumption of equal weights on the edges, there are of course many additional model extensions that could be of interest both from an applied and from a theoretical point of view. For example, the present model has no clustering, meaning that it does not contain any short cycles. Many empirical networks however have large clustering and it would therefore be relevant to try to include this in the model. Furthermore, it would be natural to allow for dependence between weights of edges at the same vertex. For example, one could think of the total weight of a vertex as being fixed and distributed randomly in some way among the edges. Also the choice of the transmission probability could be extended. Here we have frequently taken $\pi(w)=1-(1-s)^w$, where $s$ denotes the per-contact transmission probability. More generally, the per-contact transmission probability could be modeled as a random variable $S$, defined to be i.i.d.\ among different pairs of connected vertices (with for instance a Beta-distribution), and the unconditional infection probability along an edge with weight $w$ is then given by $\pi(w)=E[1-(1-S)^w]$.

Finally, there are of course other ways of constructing a weighted network with prescribed degree distribution. One suggestion is the following. First perform the standard configuration model with degree distribution $\{p(d)\}$. Then go through all edges and give weights to them independently in that an edge between two vertices with degrees $d_1$ and $d_2$ is given a weight from a distribution $q(w|d_1, d_2)$. In contrast to the model of the present paper, which uses the weights to construct the network, this model uses the network to prescribe the weights. As a consequence, this model will not have any degree correlation, which is a drawback, since many empirical networks exhibit such correlation. Indeed, social networks tend to have positive degree correlation, while some technological networks tend to have negative degree correlation; see \cite{N1,N2}. The model of the present paper only seems to be able to produce positive degree correlation (although this remains to be verified). It would be of interest to find a natural model for a weighted graph that can also produce negative degree correlation.

\section*{Acknowledgements}

All authors are grateful to Riksbankens Jubileumsfond (The Bank of Sweden Tercentenary Foundation) for financial support.

\end{document}